\documentclass[11pt]{amsart}
\newtheorem{theorem}{Theorem}[section]
\newtheorem{corollary}[theorem]{Corollary}
\newtheorem{lemma}[theorem]{Lemma}
\newtheorem{proposition}[theorem]{Proposition}

\theoremstyle{definition}
\newtheorem{definition}[theorem]{Definition}
\theoremstyle{remark}
\newtheorem{remark}[theorem]{Remark}
\numberwithin{equation}{section}

\newcommand{\CC}{\mathcal C}
\newcommand{\B}{\mathbb B}
\newcommand{\C}{\mathbb C}
\newcommand{\R}{\mathbb R}

\begin{document}

\title[Model almost complex manifolds]
{On the geometry of model almost complex manifolds with boundary}

\author{Herv\'e Gaussier and Alexandre Sukhov}

\address{\begin{tabular}{lll}
Herv\'e Gaussier & & Alexandre Sukhov\\
C.M.I. & & U.S.T.L. \\
39, rue Joliot-Curie, & & Cit\'e Scientifique \\
13453 Marseille Cedex 13 & & 59655 Villeneuve d'Ascq Cedex\\
FRANCE & & FRANCE\\
 & & \\
{\rm gaussier@cmi.univ-mrs.fr} & &
{\rm sukhov@agat.univ-lille1.fr}
\end{tabular}
}

%\email{coupet@cmi.univ-mrs.fr, \ gaussier@cmi.univ-mrs.fr,
%\ sukhov@agat.univ-lille1.fr} 


\date{\number\year-\number\month-\number\day}

\begin{abstract}
We study some special almost complex structures on strictly
pseudoconvex domains in $\R^{2n}$. They appear naturally as limits
under a nonisotropic scaling procedure and play a role of model
objects in the geometry of almost complex manifolds with boundary.  We
determine explicitely some geometric invariants of these model
structures and derive necessary and sufficient conditions for their
integrability.  As applications we prove a boundary extension and a
compactness principle for some elliptic diffeomorphisms between
relatively compact domains.
\end{abstract}

\maketitle

\section*{Introduction and main results}

The development of almost complex geometry started in the second half
of the twentieth century.  Due to the fast expansion of complex
geometry, the leadind quest, characterized by the striking theorem of
Newlander and Nirenberg~\cite{new-nir}, was to try to endow a manifold
with a complex structure. The main trespass of non integrable almost
complex manifolds was the lack of ``complex'' coordinates, essential
in both the geometric study (study of Stein manifolds,...) and the
analytic study (study of the Bergman kernel, $L^2$ estimates in
pseudoconvex domains,...).  At the same time, Nijenhuis and Woolf lead
a capital study of almost complex manifolds~\cite{ni-wo}. Their paper
may probably be considered as the starting point of the current
development of the field. Viewing almost complex maps as solutions of
non linear elliptic operators they deduced regularity results and
stability phenomena for such maps and included the almost complex
geometry in a geometric theory of elliptic partial differential operators.

In the past twenty years, symplectic geometry has been the field of
many developments. For instance, M.Gromov proved the Nonsqueezing
theorem, stating that there is no symplectic embedding of a ball into
a ``complex'' cylinder with smaller radius, and A.Floer proved
Arnold's conjecture on the number of fixed points for a symplectic
diffeomorphism in certain manifolds, developing Morse theory on
infinite-dimensional spaces. A main step in most of the recent
developments in symplectic geometry relies on the existence of
holomorphic discs. Given a symplectic form, the set of compatible
almost complex structures is a non-empty contractible oriented
manifold. As observed by M.Gromov, the space of complex curves in an
almost complex manifold tells much
information about the structure of the manifold. Symplectic invariants
of the manifold appear as invariants of the cobordism class of the
moduli space of holomorphic curves for any compatible almost complex
structure.  Underlying almost complex structures in symplectic
geometry are involved, in the issue of Nijenhuis-Woolf's
work~\cite{ni-wo}, by geometric properties of elliptic
operators. Fredholm theory provides the moduli space of holomorphic
curves or spheres with a structure of an oriented manifold, and with a
cobordism between moduli space of two distinct almost complex
structures. One views therefore almost complex manifolds as natural
manifolds for deformation theory (both of the structure and of the
associated complex curves).  The pertinence of this point of view is
dependent of some compactness principle for associated complex curves. These
compactness phenomena rely mainly on the Sobolev theory.

\vskip 0,1cm
Our paper is dedicated to the study of strictly pseudoconvex domains
in almost complex manifolds. They appear naturally in Gromov's theory.
Our approach is based on some deformation of almost complex manifolds
with boundary.  Inspired by the well-known methods of complex analysis
and geometry~\cite{pi}, we perform non isotropic dilations, naturally
associated with the geometric study of strictly convex domains in the
euclidean space.  The cluster set of deformed structures forms a
smooth non trivial manifold of {\it model} almost complex structures
on the euclidean space, containing the standard structure.  Such
nonisotropic deformations are relevant for several problems of
geometric analysis on almost complex manifolds. In the previous
paper~\cite{ga-su} we used this method to obtain lower estimates of the
Kobayashi-Royden infinitesimal metric near the boundary of a strictly
pseudoconvex domain. These estimates are one of our main technical tools
in the present paper. In the present paper we consider two distinct problems.
The first problem affects the elliptic boundary
regularity of diffeomorphisms.  In the spirit of Fefferman's
theorem~\cite{fe} on the smooth extension of biholomorphisms between smooth
strictly pseudoconvex domains, Eliashberg raised the following
question. How does a symplectic diffeomorphism of the ball effect the contact
structure of the sphere ?  One approach consists in considering a
compatible almost complex structure on the ball and study the
extension of the push forward structure under the action of the
diffeomorphism. This leads to an elliptic boundary regularity
problem. Generically, this structure does not extend up to the sphere,
since there exist symplectic diffeomorphisms which do not extend up to
the sphere.  However we prove that under some natural curvature
conditions, the extension of this structure implies the smooth
extension of the diffeomorphism up to the boundary. More precisely we
have~:

\begin{theorem}\label{theo-fefferman}
Let $D$ and $D'$ be two smooth relatively compact domains in real manifolds.
Assume that $D$ admits an almost complex structure $J$ smooth on $\bar D$ and
 such that $(D,J)$ is strictly pseudoconvex. Then a smooth
diffeomorphism $f:  D  \rightarrow  D'$ extends to a smooth
diffeomorphism between $\bar D$ and $\bar D'$  if and only
if the direct image $f_*(J)$ of $J$ under
$f$ extends smoothly on $ \bar D'$ and $(D', f_*(J))$ 
is strictly pseudoconvex.
\end{theorem}

Theorem~\ref{theo-fefferman} was proved in real dimension four in a previous
paper~\cite{co-ga-su}. In that situation, one can find a normalization
of the structure such that the cluster set for dilated structures (note
that dilations depend deeply on a choice of coordinates) is reduced to
the standard integrable structure. In the general case, the manifold
of model structures is non trivial, making the geometric study of model
structures consistent. Thus in the present paper we give a definitive
result, generalizing Fefferman's theorem (dealing with the case where
$D$ and $D'$ are equipped with the standard structure of $\C^n$).
Theorem~\ref{theo-fefferman} gives a criterion for the boundary
extension of a diffeomorphism between two smooth manifolds, under the
assumption that the source manifold admits an almost complex structure.
So it can be viewed as a geometric version of the elliptic regularity.

The second problem concerns a compactness phenomenon for some diffeomorphisms.
As this should be expected from the above general presentation, the
study of the compactness of diffeomorphisms is transformed into the study of
the compactness of induced almost complex structures, and consequently to
an elliptic problem. We prove the following compactness principle~:

\begin{theorem}\label{wr}
Let $(M,J)$ be an almost complex manifold, not equivalent to a model
domain. Let $D=\{r<0\}$ be a relatively compact domain in a smooth
manifold $N$ and let $(f^\nu)_\nu$ be a sequence of diffeomorphisms from $M$
to $D$. Assume that

$(i)$ the sequence $(J_\nu:=f^\nu_*(J))_\nu$ extends smoothly up to $\bar{D}$
and is compact in the $C^2$ convergence on $\bar{D}$,

$(ii)$ the Levi forms of $\partial D$ , $\mathcal L^{J_\nu}(\partial D)$ are
uniformly bounded from below (with respect to $\nu$) by a positive constant.

Then the sequence $(f^\nu)_\nu$ is compact in the compact-open topology on $M$.
\end{theorem}

The paper is organized as follows. In the preliminary section one, we recall some basic notions of almost complex geometry.
Section two is crucial. We introduce model almost complex structures and study
their geometric properties.
Section three contains a technical background necessary for the proof of
Theorem~\ref{theo-fefferman}. It mainly concerns properties of lifts of
almost complex structures to tangent and cotangent bundles of a manifold.
We use it to prove the boundary regularity of pseudoholomorhic discs attached
to a totally real submanifold by means of geometric bootstrap arguments.
In section four we describe nonisotropic deformations of strictly
pseudoconvex almost complex manifolds with boundary. This allows to reduce
the study of these manifolds to model structures of section one.
In section five we prove Theorem~\ref{theo-fefferman}. Our approach is
inspired by the approach of Nirenberg-Webster-Yang~\cite{nir-we-ya},
~\cite{pi-kh},~\cite{for},~\cite{tu,tu01}.
Finally in section six we prove Theorem~\ref{wr}.

\vskip 0,1cm
\noindent{\it Acknowledgments.} The authors wish to thank Bernard Coupet
and Jean-Pierre Rosay for their comments and fruitfull discussions.%\vskip 0,1cm
%The paper is organized as follows. In the first section, we recall the
%background notions. In the second section we
\section{Preliminaries}
An almost complex structure on a smooth ($\mathcal C^\infty$) real
$(2n)$-dimensional manifold $M$ is a
$\mathcal C^\infty$-field $J$ of complex linear structures on the tangent
bundle $TM$ of $ M$.
We call the pair $(M,J)$ an {\it almost complex manifold}.
We denote by $J_{st}$ the standard structure in $\R^{2n}$ and by
$\B$ the unit ball in $\R^{2n}$.
An important special case of an almost complex manifold is a bounded
domain $D$ in $\C^n$ equipped with an almost complex structure $J$,
defined in a neighborhood of $\bar{D}$, and sufficiently close to the
standard structure $J_{st}$ in the $\mathcal C^2$ norm on
$\bar{D}$ and every almost complex manifold may be represented locally
in such a form.
More precisely, we have the following Lemma. 

\begin{lemma}
\label{suplem1}
Let $(M,J)$ be an almost complex manifold. Then for every point $p \in
M$ and every $\lambda_0 > 0$ there exist a neighborhood $U$ of $p$ and a
coordinate diffeomorphism $z: U \rightarrow \mathbb B$ such that
$z(p) = 0$, $dz(p) \circ J(p) \circ dz^{-1}(0) = J_{st}$  and the
direct image $\hat J = z_*(J)$ satisfies $\vert\vert \hat J - J_{st}
\vert\vert_{\CC^2(\overline {\mathbb B})} \leq \lambda_0$.
\end{lemma}
\proof There exists a diffeomorphism $z$ from a neighborhood $U'$ of
$p \in M$ onto $\mathbb B$ satisfying $z(p) = 0$ and $dz(p) \circ J(p)
\circ dz^{-1}(0) = J_{st}$. For $\lambda > 0$ consider the dilation
$d_{\lambda}: t \mapsto \lambda^{-1}t$ in $\C^n$ and the composition
$z_{\lambda} = d_{\lambda} \circ z$. Then $\lim_{\lambda \rightarrow
0} \vert\vert (z_{\lambda})_{*}(J) - J_{st} \vert\vert_{\CC^2(\overline
{\mathbb B})} = 0$. Setting $U = z^{-1}_{\lambda}(\mathbb B)$ for
$\lambda > 0$ small enough, we obtain the desired statement. \qed

\vskip 0,1cm
Every complex one form $w$ on $M$ may be uniquely
decomposed as $w=w_{(1,0)} + w_{(0,1)}$, where $w_{(1,0)} \in
T^*_{(1,0)}M$ and $w_{(0,1)} \in T^*_{(0,1)}M$, with respect to the structure
$J$.   
This enables to define the operators $\partial_J$ and
$\bar{\partial}_J$ on the space of smooth functions defined on
$M$~: given a complex smooth function $u$ on $M$, we set $\partial_J u =
du_{(1,0)}$ and $\bar{\partial}_Ju = du_{(0,1)}$.

\subsection{Real submanifolds in an almost complex manifold}
Let $\Gamma$ be a real smooth submanifold in $M$.
We denote by $H^J(\Gamma)$ the $J$-holomorphic tangent bundle
$T\Gamma \cap JT\Gamma$. Then $\Gamma$ is {\it totally real}
if $H^J(\Gamma) = \{0\}$ and $\Gamma$ is {\it $J$-complex} if
$T\Gamma = H^J(\Gamma)$.

If $\Gamma$ is a real
hypersurface in $M$ defined by $\Gamma=\{r=0\}$ and $p \in \Gamma$
then by definition 
$$
H_p^J(\Gamma) = \{v \in T_pM : dr(p)(v) =
dr(p)(J(p)v) = 0\} = \{v \in T_pM :
\partial_Jr(p) (v-iJ(p)v) = 0\}.
$$
We recall the notions of the Levi form~:

\begin{definition}\label{DEF}
Let $\Gamma=\{r=0\}$ be a smooth real hypersurface in $M$ 
($r$ is any smooth defining function of $\Gamma$) and let $p \in \Gamma$. 

$(i)$ The {\sl Levi form} of $\Gamma$ at $p$ is the map defined on
$H^J_p(\Gamma)$  by $\mathcal L^J(\Gamma)(X_p) = J^\star dr[X,JX]_p$,
where the vector field $X$ is any section of the $J$-holomorphic tangent
bundle  $H^J \Gamma$ such that $X(p) = X_p$.

$(ii)$ A real smooth hypersurface $\Gamma=\{r=0\}$ in $M$ is 
{\sl strictly $J$-pseudoconvex} if its Levi form $\mathcal L^J(\Gamma)$ 
is positive definite on $H^J(\Gamma)$.

$(iii)$ If $r$ is a $\CC^2$ function on $M$ then the Levi 
form of $r$ is defined on $TM$ by $\mathcal L^J(r)(X):=-d(J^\star dr)(X,JX)$.

$(iv)$ A $\CC^2$ real valued function $r$ on $M$ is
$J$-plurisubharmonic on $M$ (resp. strictly $J$-plurisubharmonic)
if and only if $\mathcal L^J(r)(X) \geq 0$ for every $X \in TM$ (resp.
$\mathcal L^J(r)(X) > 0$ for every $X \in TM \backslash \{0\}$).
\end{definition}

\subsection{Local representation of holomorphic discs}
A smooth map $f$ between two almost complex manifolds $(M',J')$ and
$(M,J)$ is holomorphic if its differential satisfies the following
holomorphy condition~: $df \circ J = J' \circ df$ on $TM$. In case
$(M',J') = (\Delta, J_{st})$ the map $f$ is called a $J$-holomorphic
disc. We denote by $\zeta$ the complex variable in $\C$.  In view of
Lemma~\ref{suplem1}, the holomorphy condition is usually written as
$$ \frac{\partial f}{\partial \bar{\zeta}} + Q_J(f) \frac{\partial
f}{\partial \zeta} = 0,
$$ where $Q=(J_{st}+J)^{-1}(J_{st}-J)$ (see~\cite{si}).
However, in view of Lemma~\ref{suplem1} a basis $w:=(w^1,\dots,w^n)$
of $(1,0)$ forms on $M$ may be locally written as $w^j=dz^j +
\sum_{k=1}^nA_{j,k}(z,\bar{z})d\bar{z}$ where $A_{j,k}$ is a smooth
function.  The disc $f$ being $J$-holomorphic if $f^*(w^j)$ is a
$(1,0)$ form for $j=1,\dots,n$ (see~\cite{ch1}), $f$ satisfies the
following equation on $\Delta$~:
\begin{equation}\label{eq1}
\frac{\partial f}{\partial \bar{\zeta}} + A(f)
\overline{\frac{\partial f}{\partial \zeta}} = 0,
\end{equation}
where $A=(A_{j,k})_{1 \leq j,k \leq n}$. 
We will use this second equation to characterize the $J$-holomorphy
in the paper.

\section{Model almost complex structures}

The scaling process in complex manifolds deals with deformations of
domains under holomorphic transformations called dilations. The usual
nonisotropic dilations in complex manifolds, associated with strictly
pseudoconvex domains, provide the unit ball (after biholomorphism) as the
limit domain. In almost complex manifolds dilations are generically no more
holomorphic with respect to the ambiant structure. The scaling process
consists in deforming both the structure and the domain.
This provides, as limits, a quadratic domain and a linear deformation
of the standard structure in $\R^{2n}$, called {\it model structure}.
We study some invariants of such structures.
Let $(x^1,y^1,\dots,x^n,y^n)=(z^1,\dots,z^n)=('z,z^n)$ denote the canonical
coordinates of $\R^{2n}$.
\begin{definition}\label{def-model}
{\it Let $J$ be an almost complex structure on $\C^n$. We call
$J$ a {\rm model structure} if $J(z) = J_{st} + L(z)$ where $L$ is given by
a linear matrix $L=(L_{j,k})_{1 \leq j,k \leq n}$ such that $L_{j,k} = 0$ for
$1 \leq j \leq 2n-2, \ 1 \leq k \leq 2n$, $L_{j,k}=0$ for $j,k = 2n-1,2n$  and
$L_{n,k} = \sum_{l=1}^{n-1} (a_l^k z^l + \bar{a}_l^k \bar{z}^l)$,
$a_l^k \in \C$.}
\end{definition}
The complexification $J_\C$ of a model structure $J$ can be written
as a $(2n \times 2n)$ complex matrix
\begin{equation}\label{complex}
J_\C=\left(
\begin{array}{ccccccc}
i      & 0      & 0      & 0      & \cdots & 0      & 0 \\
0      & -i     & 0      & 0      & \cdots & 0      & 0 \\
0      & 0      & i      & 0      & \cdots & 0      & 0 \\
0      & 0      & 0      & -i     & \cdots & 0      & 0 \\
\cdots & \cdots & \cdots & \cdots & \cdots & \cdots & \cdots \\
0      & \tilde{L}_{2n-1,1} & 0 & \tilde{L}_{2n-1,2} & \cdots & i & 0\\
\tilde{L}_{2n,1} & 0 & \tilde{L}_{2n,3} & 0 & \cdots &  0 & -i
\end{array}
\right),
\end{equation}
where $\tilde{L}_{2n-1,k}(z,\bar{z}) = \sum_{l=1,\ l \neq k}^{n-1}
(\alpha_l^{k} z^l + \beta_l^{k} \bar{z}^l)$
with $\alpha_l^k,\ \beta_l^k \in \C$.
Moreover, $\tilde{L}_{2n,2k-1} = \overline{\tilde{L}_{2n-1,2k}}$.
%The almost complex condition $J^2 = -I$ for a model structure implies that
%$L_{,j} \circ J_{st} + J_{st} \circ L_{n,j} = 0$ for $j=1,\dots,n-1$.
%In particular if $v$ is a complex vector in $\C$ then
%$L_{n,j}(v) = \sum_{l=1}^{n-1} (\alpha_l^j z^l + \beta_l^j \bar{z}_l)\bar v$.

\vskip 0,1cm
With a model structure we associate model domains.
%\begin{definition}
%Let $J$ be an almost complex structure
%in $\C^n$ and let $\mathbb H_P=\{z \in \C^n :
%z^n + \bar{z}_n + P('z,'\bar{z}) < 0\}$ where $P$ is a homogeneous polynomial
%of degree two. We call $(\mathbb H_P,J)$ a {\it model domain} if $P$
%is $J$-pluri-subharmonic on $\C^{n-1}$ and $P$ is strictly
%$J$-pluri-subharmonic in a neighborhood of the origin.
%\end{definition}

\begin{definition}
Let $J$ be a model structure on $\C^n$ and $D=\{z \in \C^n : Re z^n +
P_2('z,'\bar{z})<0\}$, where $P_2$ is homogeneous second degree real
polynomial on $\C^{n-1}$.  The pair $(D,J)$ is called a {\it model
domain} if $D$ is strictly $J$-pseudoconvex in a neighborhood of the
origin.
\end{definition}

The aim of this Section is to define the complex hypersurfaces for model
structures in $\R^{2n}$.

Let $J$ be a model structure on $\R^{2n}$ and let $N$ be a germ of a
$J$-complex hypersurface in $\R^{2n}$.
\begin{proposition}\label{prop-hyp}
\hfill

$(i)$ The model structure $J$ is integrable if and only if
$ \tilde{L}_{2n-1,j}$ satisfies the compatibility conditions
$$
\frac{\partial \tilde{L}_{2n-1,k}}{\partial \bar{z}^j} =
\frac{\partial \tilde{L}_{2n-1,j}}{\partial \bar{z}^k}
$$
for every $1 \leq j,k \leq n-1$. 

In that case there exists a global diffeomorphism of $\R^{2n}$
which is $(J,J_{st})$ holomorphic. In that case the germs of any $J$-complex
hypersurface are given by one of the two following forms~:

\hskip 0,5cm $(a)$ $N= A \times \C$ where $A$ is a germ of a $J_{st}$-complex
hypersurface in $\C^{n-1}$,

\hskip 0,5cm $(b)$ $N =\{('z,z^n) \in \C^n : z^n = \frac{i}{4}
\sum_{j=1}^{n-1}\bar{z}^j\tilde{L}_{2n-1,j}('z,'\bar{z})
+\frac{i}{4}
\sum_{j=1}^{n-1}\bar{z}^j \tilde{L}_{2n-1,j}('z,0))
+ \tilde{\varphi}('z)\}$
where $\tilde{\varphi}$ is a holomorphic function locally defined in
$\C^{n-1}$.

\vskip 0,1cm
$(ii)$ If $J$ is not integrable then $N= A \times \C$ where $A$ is a germ
of a $J_{st}$-complex hypersurface in $\C^{n-1}$.
\end{proposition}
\vskip 0,2cm
\noindent{\bf Proof of Proposition~\ref{prop-hyp}}. Let $N$ be a germ
of a $J$-complex hypersurface in $\R^{2n}$.
If $\pi:\R^{2n} \rightarrow \R^{2n-2}$ is the projection on the $(2n-2)$ first
variables, it follows from Definition~\ref{def-model}, or similarly
from condition~(\ref{complex}) that
$\pi(T_zN)$ is a $J_{st}$-complex hypersurface in $\C^{n-1}$.
%Indeed, for every $v=(v^1,w^1,\dots, v^{n-1}, w^{n-1}, v^n,w^n) \in T_zN$,
%$\pi(J_z(v)) = (-w^1,v^1,\dots, -w^{n-1},v^{n-1})$.

It follows that either $dim_\C\pi(N) = n-1$ or $dim_\C\pi(N) = n-2$. 
\vskip 0,1cm
\noindent{\it Case one : $dim_\C\pi(N) = n-1$}.
We prove the following Lemma~:

\begin{lemma}\label{lem-hyp}
There is a local holomorphic function $\tilde{\varphi}$ in $\C^{n-1}$ such that
$N =\{('z,z^n) : z^n = \frac{i}{4}
\sum_{j=1}^{n-1}\bar{z}^j\tilde{L}_{2n-1,j}('z,'\bar{z})
+\frac{i}{4}
\sum_{j=1}^{n-1}\bar{z}^j \tilde{L}_{2n-1,j}('z,0))
+ \tilde{\varphi}('z)\}.$
\end{lemma}

{\it Proof of Lemma~\ref{lem-hyp}}. A germ $N$ can be represented
as a graph $N=\{z^n = \varphi('z,'\bar{z})\}$ where $\varphi$ is a smooth local
complex function. Hence
$T_zN=\{v_n = \sum_{j=1}^{n-1}(\frac{\partial \varphi}{\partial z^j}('z)v_j 
+ \frac{\partial \varphi}{\partial {\bar z}^j}('z) \bar{v}_j)\}$.
A vector $v=(x^1,y^1,\dots,x^n,y^n)$ belongs to $T_zN$ if and only if
the complex components $v^1:=x^1+iy^1,\dots,v^n:=x^n + i y^n$ satisfy
\begin{equation}\label{tan}
iv_n = i\sum_{j=1}^{n-1}(\frac{\partial \varphi}{\partial z^j}('z)v_j +
\frac{\partial \varphi}{\partial \bar{z}^j}('z)\bar{v}_j).
\end{equation}
Similarly, the vector $J_zv$ belongs to $T_zN$ if and only if
\begin{equation}\label{Jtan}
\sum_{j=1}^{n-1}\tilde{L}_{2n,2j-1}('z) \bar{v}_j + i v_n =
i(\sum_{j=1}^{n-1} \frac{\partial \varphi}{\partial z^j}('z)
v_j -\sum_{j=1}^{n-1}\frac{\partial \varphi}{\partial \bar{z}^j}('z)\bar{v}_j).
\end{equation}
It follows from (\ref{tan}) and (\ref{Jtan}) that $N$ is $J$-complex if and
only if
$$
\sum_{j=1}^{n-1}(\tilde{L}_{2n,2j-1}('z)
\bar{v}_j + 2i \frac{\partial \varphi}{\partial \bar{z}^j}(z)\bar{v}_j) = 0
$$
for every $'v \in \C^{n-1}$, or equivalently if and only if
$$
\tilde{L}_{2n,2j-1} = -2i \frac{\partial \varphi}{\partial \bar{z}^j}
$$
for every $j=1,\cdots,n-1$.
This last condition is equivalent to the compatibility
conditions
\begin{equation}\label{compat}
\frac{\partial \tilde{L}_{2n,2j-1}}{\partial \bar{z}^k} =
\frac{\partial \tilde{L}_{2n,2k-1}}
{\partial \bar{z}^j}\ {\rm for}\ j,k = 1,\cdots,n-1.
\end{equation}
In that case
there exists a local holomorphic function $\tilde{\varphi}$ in $\C^{n-1}$
such that
$$
\varphi('z,'\bar{z})=\frac{i}{2}
\sum_{j=1}^{n-1}\bar{z}^j(\sum_{k \neq j}\alpha_k^j z^k)
-\frac{i}{2}
\sum_{j=1}^{n-2}\bar{z}^j(\sum_{k > j}\beta_k^j \bar{z}^k)
+ \tilde{\varphi}('z),
$$
meaning that such $J$-complex hypersurfaces are parametrized by holomorphic
functions in the variables $'z$.
Moreover we can rewrite $\varphi$ as
$$
\varphi('z,'\bar{z})=\frac{i}{4}
\sum_{j=1}^{n-1}\bar{z}^j\tilde{L}_{2n-1,j}('z,'\bar{z})
+\frac{i}{4}
\sum_{j=1}^{n-1}\bar{z}^j \tilde{L}_{2n-1,j}('z,0))
+ \tilde{\varphi}('z).
$$ \qed

%We come back to the proof of Proposition~\ref{prop-hyp}.
We also have the following
\begin{lemma}\label{lem-form}
The $(1,0)$ forms of $J$ have the form
$\omega = \sum_{k=1}^n c_kdz^k -\frac{i}{2}c_n\sum_{k=1}^{n-1}
\tilde{L}_{2n-1,k} d\bar{z}^k$ with complex numbers $c_1,\dots,c_n$.
\end{lemma}

{\it Proof of Lemma~\ref{lem-form}}.
Let $X=\sum_{k=1}^n(x_k \frac{\partial}{\partial z^k} +
y_k\frac{\partial}{\partial \bar{z}^k})$ be a $(0,1)$ vector field.
In view of (\ref{complex}), we have~:
$$
J_\C(X) = -iX \Leftrightarrow \left\{
\begin{array}{lll}
x_k & = & 0, \ \ {\rm for k =1,\dots,n-1}\\
 & & \\
x_n & = & \frac{i}{2} \sum_{k=1}^{n-1}y_k \tilde{L}_{2n-1,k}.
\end{array}
\right.
$$
Hence the $(0,1)$ vector fields are given by
$$
X=\sum_{k=1}^ny_kd\bar{z}^k +
\frac{i}{2}dz^n \sum_{k=1}^{n-1}y_k\tilde{L}_{2n-1,k}.
$$
A $(1,0)$ form $\omega=\sum_{k=1}^n(c_kdz^k + d_nd\bar{z}^k)$ 
satisfying $\omega(X)=0$ for every $(0,1)$ vector field $X$ it satisfies
$d_n = 0$ and $d_k + (i/2)c_n \tilde{L}_{2n-1,k}=0$ for every
$k=1,\dots,n-1$.
This gives the desired form for the $(1,0)$ forms on $\C^n$. \qed
\vskip 0,1cm
Consider now the global diffeomorphism of $\C^n$ defined by
$$
F('z,z^n) = ('z,z^n - \frac{i}{4}
\sum_{j=1}^{n-1}\bar{z}^j\tilde{L}_{2n-1,j}('z,'\bar{z})
-\frac{i}{4}
\sum_{j=1}^{n-1}\bar{z}^j \tilde{L}_{2n-1,j}('z,0)).
$$
The map $F$ is $(J,J_{st})$ holomorphic if and only if $F^*(dz^k)$
is a $(1,0)$ form with respect to $J$, for every $k=1,\dots,n$.

Then $F^*(dz^k) = dz^k$ for $k=1,\dots,n-1$ and
$$
\begin{array}{lll} 
F^*(dz^n) & = & \displaystyle
dz^n + \sum_{k=1}^{n-1} \frac{\partial F_n}{\partial z^k} dz^k
+ \sum_{k=1}^{n-1} \frac{\partial F_n}{\partial \bar{z}^k} d\bar{z}^k\\
& & \\
 & = & \displaystyle
dz^n + \sum_{k=1}^{n-1} \frac{\partial F_n}{\partial z^k} dz^k\\
& & \ \ \ \ \ - \displaystyle \frac{i}{4} \sum_{k=1}^{n-1}
(\tilde{L}_{2n-1,k}('z,'\bar{z})+
\sum_{j \neq k}\bar{z}^j
\frac{\partial \tilde{L}_{2n-1,j}}{\partial \bar{z}^k}('z,'\bar{z})
+ \tilde{L}_{2n-1,k}('z,'0)) d\bar{z}^k.
\end{array}$$
By the compatibility condition~(\ref{compat}) we have
$$
\begin{array}{lll}
F^*(dz^n) & = & \displaystyle
dz^n + \sum_{k=1}^{n-1} \frac{\partial F_n}{\partial z^k} dz^k
-\frac{i}{4} \sum_{k=1}^{n-1}
(\tilde{L}_{2n-1,k}('z,'\bar{z}) +
\tilde{L}_{2n-1,k}('0,'\bar{z}) +
\tilde{L}_{2n-1,k}('z,'0)) d\bar{z}^k\\
& & \\
 & = & \displaystyle dz^n 
- \frac{i}{2} \sum_{k=1}^{n-1}
\tilde{L}_{2n-1,k}('z,'\bar{z}) d\bar{z}^k
+ \sum_{k=1}^{n-1} \frac{\partial F_n}{\partial z^k} dz^k.
\end{array}
$$
These equalities mean that $F$ is a local $(J,J_{st})$-biholomorphism of
$\C^n$, and so that $J$ is integrable.

\vskip 0,1cm
\noindent{\it Case two : $dim_\C\pi(N) = n-2$}. In that case we may write
$N=\pi(N) \times \C$, meaning that $J$-complex hypersurfaces
are parametrized by $J_{st}$-complex hypersurfaces of $\C^{n-1}$.

\vskip 0,1cm
We can conclude now the proof of Proposition~\ref{prop-hyp}. We proved in
Case one that if there exists a $J$-complex hypersurface in $\C^n$ such
that $dim \pi(N) = n-1$ (this is equivalent to the compatibility conditions
(\ref{compat})) then $J$ is integrable.
Conversely, it is immediate that if $J$
is integrable then there exists a $J$-complex hypersurface whose form
is given by Lemma~\ref{lem-hyp} and hence that the compatibility conditions
(\ref{compat}) are satisfied. This gives part $(i)$ of
Proposition~\ref{prop-hyp}.

To prove part $(ii)$, we note that if $J$ is not integrable then in view
of part $(i)$ the form of any $J$-complex hypersurface is given by Case two.
 \qed

\section{Almost complex structures and totally real submanifolds}
We recall that if $(M,J)$ is an almost complex manifold then a submanifold $N$
of $M$ is {\it totally real} if $TN \cap J(TN) = \{0\}$. 

\subsection{Conormal bundle of a submanifold in $(M,J)$}
The conormal bundle of a strictly $J$-pseudoconvex hypersurface in $M$ provides
 an important example of a totally real submanifold in the cotangent bundle
$T^*M$.
%Let $(M,J)$ and $(M',J')$  be almost complex manifolds. We denote by
%$T^*M$ the cotangent bundle of $M$.
More precisely, following Sato~
(see~\cite{ya-is}) let $\widehat J$ denote the
complete lift of $J$. If $J$ has
components $J_i^{\, h}$ then, in the matrix form we have

$$
\widehat{J} = \left(
\begin{array}{cll}
J_i^{\, h} & & 0 \\
 p_a \left ( \frac{\partial J_i^{\, a}}{\partial x^j} -
\frac{\partial J_j^{\, a}}{\partial x^i} \right ) & & J_h^{\, i}
\end{array}
\right),
$$
relative to local canonical coordinates $(x,p)$ on $T^*M$.

Let $N$ denote the Nijenhuis tensor and let $\gamma(NJ)$
be the $(1,1)$ tensor on $T^*M$ defined in coordinates $(x,p)$ by~:
$$
\gamma (NJ) = \left(
\begin{array}{cll}
0 & & 0\\
 p_a (NJ)_{ji}^{\,\,a} & & 0 
\end{array}
\right).
$$
Then the $(1,1)$ tensor defined on
$T^*M$ by $\tilde{J} := \widehat J
+(1/2)\gamma(NJ)$ defines an almost complex structure on the cotangent bundle
$T^*M$. 
Moreover, if $f$ is a biholomorphism
between $(M,J)$ and $(M',J')$ then the cotangent map
$\tilde f:= (f,{}^tdf^{-1})$ is a biholomorphism between $(T^*M,\tilde J)$ and
$(T^*M',\tilde J')$.
%The main properties of $\tilde J$ are given in
%the following proposition (see \cite{co-ga-su}, Proposition~6.2)~:

%\begin{proposition}\label{prop-lift}
%$(i)$ If $f$ is a biholomorphism between $(M,J)$ and $(M',J')$
%then the cotangent map $\tilde f:= (f,{}^tdf^{-1})$ is a biholomorphism
%between $(T^*M,\tilde J)$ and $(T^*M',\tilde J')$.

%$(ii)$ If $(J_\varepsilon)_\varepsilon$ is a small deformation of the
%standard structure on $\C^n$ then
%$\tilde{J}_\varepsilon \rightarrow J_{st}$ as $\varepsilon \rightarrow 0$
%in the $\CC^k$-norm on $T^*\C^n$ (for any $k$).
%\end{proposition}

%Let $i: T^*(M) \rightarrow  T^*_{(1,0)}(M,J) $ be the
%canonical identification between the cotangent bundle of $M$ and the space
%of $(1,0)$-forms on $TM$ with respect to $J$. 
%%Let $D$ be a smoothly relatively compact
%%domain in $M$ with boundary $\Gamma$.
If $\Gamma$ is a real submanifold in $M$, the conormal bundle
$\Sigma_J(\Gamma)$ of $\Gamma$ is the real subbundle of $T^*_{(1,0)}M$
defined by $\Sigma_J(\Gamma) = \{ \phi \in T^*_{(1,0)}M: Re \,\phi
\vert_{T\Gamma} = 0\}$.  One can identify the
conormal bundle $\Sigma_J(\Gamma)$ of $\Gamma$ with any of the
following subbundles of $T^*M$~: $N_1(\Gamma)=\{\varphi \in T^*M :
\varphi_{|T\Gamma}=0\}$ and $N_2(\Gamma)=\{\varphi \in T^*M :
\varphi_{|JT\Gamma}=0\}$.
%If $f:(D,J) \rightarrow (D',J')$
%is a biholomorphic map $\CC^1$-smooth up to $\partial D$, then its
%cotangent map $\tilde f$ is continuous
%up to $\Sigma_J(\partial D)$ and $\tilde{f}(\Sigma_J(\partial
%D)) = \Sigma_{J'}(\partial D')$.

\begin{proposition}\label{prop-tot-real}
Let $\Gamma$ be a $\CC^2$ real hypersurface in $(M,J)$.  If the Levi
form of $\Gamma$ is nondegenerate, then the bundles $N_1(\Gamma)$ and
$N_2(\Gamma)$ (except the zero section) are totally real submanifolds
of dimension $2n$ in $T^*M$ equipped with $\tilde{J}$.
\end{proposition} 
There is in fact an equivalence between the nondegeneracy of the Levi form
of $\Gamma$ and the fact that the conormal bundle of $\Gamma$ is totally
real. However we focus on the implication suitable for our purpose.
Proposition~\ref{prop-tot-real} is due to A.Tumanov~\cite{tu01} in the
integrable case. The question wether a similar result was true in the almost complex case was
asked by the second author to A.Spiro who gave a positive answer \cite{sp}. 
For completeness we give here  an alternative proof of the same fact. 

\vskip 0,1cm
\noindent{\bf Proof of Proposition~\ref{prop-tot-real}}. Let $x_0 \in \Gamma$.
We consider local coordinates $(x,p)$ for the real cotangent bundle
$T^*M$ of $M$ in a neighborhood of $x_0$.
The fiber of $N_2(\Gamma)$ is given by $c(x)J^*d\rho(x)$, where $c$ is a real
nonvanishing function. In what follows we denote $-J^*d\rho$
by $d^c_J\rho$. For every $\varphi \in N_2(\Gamma)$
we have $\varphi|_{J(T\Gamma)} \equiv 0$. It is equivalenty to prove that
$N_1(\Gamma)$ is totally real in $(T^*M,\tilde J)$ or that $N_2(\Gamma)$ is
totally real in $(T^*M,\tilde J)$. We recall
that if $\Theta=p_idx^i$ in local coordinates then $d\Theta$ defines the
canonical symplectic form on $T^*M$.
If $V,W \in T(N_2(\Gamma))$ then
$d\Theta(V,W)=0$. Indeed the projection $pr_1(V)$ of $V$ (resp. $W$) on $M$ is
in $J(T\Gamma)$ and the projection of $V$ (resp. $W$) on the fiber
annihilates $J(T\Gamma)$ by definition. It follows that $N_2(\Gamma)$ is a
Lagrangian submanifold of $T^*M$ for this symplectic form. 

Let $V$ be a vector field in $T(N_2(\Gamma)) \cap \tilde JT(N_2(\Gamma))$.
We wish to prove that $V=0$.
According to what preceeds we have $d\Theta(V,W) = d\Theta(JV,W)=0$ for every
$W \in T(N_2(\Gamma))$. We restrict to $W$ such that $pr_1(W) \in T\Gamma
\cap J(T\Gamma)$. Since $\Theta$ is defined over $x_0 \in \Gamma$ by
$\Theta = cd^c_J \rho$, then
$d\Theta = dc \wedge d^c_J \rho + cdd^c_J\rho$.
Since $d^c_J \rho(pr_1(V)) = d^c_J \rho(Jpr_1(V)) = d^c_J \rho(pr_1(V)) =
d^c_J \rho(Jpr_1(V)) =0$ it follows that
$dd^c_J \rho (pr_1(V),pr_1(\tilde{J}W)) =0$. However, by the definition of
$\tilde{J}$, we know that $pr_1(\tilde{J}W) = J pr_1(W)$.
Hence, choosing $W = V$, we obtain that $dd^c_J \rho (pr_1(V),J pr_1(V)) = 0$.
Since $\Gamma$ is strictly $J$-pseudoconvex, it follows that $pr_1(V) = 0$.
In particular, $V$ is given in local coordinates by $V=(0,pr_2(V))$.
It follows now from the form of $\tilde J$ that $JV=(0,J pr_2(V))$ (we consider
$pr_2(V)$ as a vector in $\mathbb R^{2n}$ and $J$ defined on $\R^{2n}$).
Since $N_2(\Gamma)$ is a real bundle of rank one, then $pr_2(V)$ is
equal to zero. \qed
 
\subsection{Boundary regularity of $J$-holomorphic discs}
Many geometric questions in complex analysis or in CR geometry reduce to
the study of the properties of holomorphic discs. Among these the boundary
regularity of holomorphic discs attached to a totally real submanifold
appeared as one of the essential tools in the understanding of extension
phenomena. In the almost complex setting, this is stated by
H.Hofer \cite{ho} (refered to a bootstrap argument), and
a weaker regularity is proved by E.Chirka \cite{ch1}
and by S.Ivashkovich-V.Shevchischin \cite{iv-sh}.
This can be formulated as follows~:
\begin{proposition}\label{reflection}
Let $N$ be a smooth $\mathcal C^\infty$ totally real submanifold in $(M,J)$
and let $\varphi : \Delta^+ \rightarrow M$ be $J$-holomorphic.
Assume that the cluster set of $\varphi$ on the real interval $]-1,1[$ is
contained in $N$. Then $\varphi$ is of class $\mathcal C^\infty$ on
$\Delta^+ \cup ]-1,1[$.
\end{proposition}

Here $\Delta$ denotes the unit disc in $\mathbb C$ and
$\Delta^+:=\{\zeta \in \Delta : Im(\zeta) >0\}$.

In case $N$ has a weaker regularity then the exact regularity of $\varphi$,
related to that of $N$, can be derived directly from the following proof of
Proposition~\ref{reflection}.

\vskip 0,1cm
\noindent{\bf Proof of Proposition~\ref{reflection}.} We proceed in three
steps, using a geometric bootstrap argument.

\vskip 0,1cm
\noindent{\it Step one : $1/2$-H\"older continuity}.  Since $N$ is
totally real, using a partition of unity, we may represent $N$ as the
zero set of the positive, smooth, strictly $J$-plurisubharmonic
function $\rho$ (see the details in ~\cite{co-ga-su}).

As usual we denote by $C(\varphi,]-1,1[)$ the cluster set of $\varphi$ on
$]-1,1[$; this consists of points $p \in M$ such that $p=\lim_{k \rightarrow
\infty}\varphi(\zeta_k)$ for a sequence $(\zeta_k)_k$ in $\Delta^+$ converging
to a point in $]-1,1[$. The $1/2$-H\"older extension of $\varphi$ to
$\Delta^+ \cup ]-1,1[$ is contained in the following proposition
(see Proposition 4.1 in ~\cite{co-ga-su} for its proof).

\begin{proposition}
\label{regth1}
Let $G$ be a relatively compact domain in $(M,J)$
and let $\rho$ be a strictly $J$-plurisubharmonic
function of class $\CC^2$ on $\bar{G}$.
Let $\varphi:\Delta^+ \rightarrow G$ be a
$J$-holomorphic disc such that $\rho \circ \varphi \geq 0$ on $\Delta^+$.
Suppose that $C(g,]-1,1[)$ is contained in the zero set of $\rho$.
Then $\varphi$ extends as a H\"older 1/2-continuous map on $\Delta^+ \cup
]-1,1[$.
\end{proposition}

\vskip 0,1cm
\noindent{\it Step two : The disc $\varphi$ is of class $\mathcal C^{1+1/2}$.}
 The following construction of the reflection principle
for pseudoholomorphic discs is due to Chirka \cite{ch1}. For reader's
convenience we give the details.
Let $a\in ]-1,1[$. Our consideration being local at $a$, we may assume that
$N=\R^n \subset \C^n$, $a=0$ and $J$ is a smooth almost complex structure
defined in the unit ball $\B_n$ in $\C^n$.

After a complex linear change of coordinates we may assume that
$J = J_{st} + O(\vert z \vert)$ and $N$ is given by $x + ih(x)$ where
$x \in \R^n$ and $dh(0) = 0$. If $\Phi$ is the local diffeomorphism
$x \mapsto x$, $y \mapsto y - h(x)$ then $\Phi(N) = \R^n$ and the direct
image of $J$ by $\Phi$, still denoted by $J$, keeps the form $J_{st} +
O(\vert z \vert)$. Then $J$ has a basis of $(1,0)$-forms given in the
coordinates $z$ by $dz^j + \sum_k a_{jk}d\bar z^k$; using the
matrix notation we write it in the form $\omega = dz + A(z)d\bar z$ where
the matrix function $A(z)$ vanishes at the origin. Writing
$\omega = (I + A)dx + i(I - A)dy$ where $I$ denotes the identity
matrix, we can take as a basis of $(1,0)$ forms~: $\omega' = dx +
i(I +A)^{-1}(I - A)dy = dx + iBdy$. Here the matrix function $B$ satisfies
$B(0) = I$. Since $B$ is smooth, its restriction $B_{\vert \R^n}$ on $\R^n$
admits a smooth extension $\hat B$ on the unit ball such that
$\hat B - B_{\vert \R^n} = O(\vert y \vert^k)$ for any positive integer $k$.
Consider the diffeomorphism $z^* = x + i\hat B(z) y$.
In the $z^*$-coordinates the submanifold $N$ still coincides with $\R^n$
and $\omega' = dx + iBdy = dz^* + i(B - \hat B)dy - i(d\hat B)y = dz^* +
\alpha$, where the coefficients of the form $\alpha$ vanish with infinite
order on $\R^n$. Therefore there is a basis of $(1,0)$-forms
(with respect to the image of $J$ under the coordinate diffeomorphism
$z \mapsto z^*$) of the form $dz^* + A(z^*)d\bar z^*$,
where $A$ vanishes to first order on $\R^n$ and
$\| A \|_{\mathcal C^1(\bar{\B}_n)} < < 1$.

Consider the continuous map $\psi$ defined on $\Delta$ by
$$
\left\{
\begin{array}{cccc}
\psi &=& \varphi &{\rm on}\ \Delta^+\\
& & & \\
\psi(\zeta) &=& \overline{\varphi(\bar{\zeta})} &{\rm for}\ \zeta \in \Delta^-.
\end{array}
\right.
$$

Since, in view of (\ref{eq1}) the map $\varphi$ satisfies
\begin{equation}\label{holo}
\bar \partial \varphi + A(\varphi)\overline{\partial \varphi} = 0
\end{equation}
on $\Delta^+$, the map $\psi$ satisfies the equation
$$
\bar\partial\psi(\zeta) + 
\overline{A(\varphi(\bar\zeta))}\
\overline{\partial\psi(\zeta)} = 0
$$
for $\zeta \in \Delta^-$.

Hence $\psi$ is a solution on $\Delta$ of the elliptic equation
\begin{equation}\label{elliptic}
\bar \partial \psi + \lambda(\cdot)\overline{\partial \psi} = 0
\end{equation}
where $\lambda$ is defined by $\lambda(\zeta) =
A(\varphi(\zeta))$ for
$\zeta \in \Delta^+ \cup ]-1,1[$ and $\lambda(\zeta) =
\overline{A(\varphi(\bar\zeta))}$ for
$\zeta \in \Delta^-$.
According to Step
one, the map $\lambda$ is H\"older $1/2$ continuous on $\Delta$
and vanishes on $]-1,1[$.
This implies that $\psi$ is of class $\mathcal C^{1+1/2}$ on $\Delta$
by equation~(\ref{elliptic}) (see~\cite{si,ve}).

\vskip 0,1cm
\noindent{\it Step three : geometric bootstrap.} Let $v=(1,0)$ in $\R^2$ and
consider the disc $\varphi^c$ defined on $\Delta^+$ by
$$
\varphi^c(\zeta) = (\varphi(\zeta),d\varphi(\zeta)(v)).
$$
We endow the tangent bundle $TM$ with the complete lift $J^c$ of $J$
(see~\cite{ya-is} for its definition). We recall that $J^c$ is an
almost complex structure on $TM$. Moreover, if $\nabla$ is any $J$
complex connection on $M$ (ie $\nabla J=0$) and $\bar{\nabla}$ is the
connection defined on $M$ by $\bar\nabla_XY = \nabla_YX +[X,Y]$ then
$J^c$ is the horizontal lift of $J$ with respect to
$\bar\nabla$. Another definition of $J^c$ is given in \cite{le-sz}
where this is characterized by a deformation property.
The equality between the two definitions given in \cite{ya-is} and in
\cite{le-sz} is obtained by their (equal) expression in the local
canonical coordinates on $TM$~:
$$
J^c=\left(
\begin{array}{ccc}
J_i^h & & (0)\\
  & &    \\
t^a\partial_a J_i^h & & J_i^h
\end{array}
\right).
$$
(Here $t^a$ are fibers coordinates).

\begin{lemma}\label{tot-real}
$TN$ is a totally real submanifold in $(TM,J^c)$.
\end{lemma}

\noindent{\bf Proof of Lemma~\ref{tot-real}.} 
Let $X \in T(TN) \cap J^c(T(TN))$. If $X=(u,v)$ in the trivialisation
$T(TM) = TM \oplus TM$ then $u \in TN \cap J(TN)$,
implying that $u=0$. Hence $v  \in TN \cap J(TN)$,
implying that $v=0$. Finally, $X=0$. \qed

\vskip 0,1cm
The cluster set $C(\varphi^c,]-1,1[)$ is contained in the smooth
submanifold $TN$ of $TM$. 
Applying Step two to $\varphi^c$ and $TN$ we prove that the first derivative
of $\varphi$ with respect to $x$ ($x+iy$ are the standard coordinates on $\C$)
is of class $\mathcal C^{1+1/2}$ on $\Delta^+ \cup ]-1,1[$. 
The $J$-holomorphy equation~(\ref{holo}) may be written as
$$
\frac{\partial \varphi}{\partial y} = J(\varphi)
\frac{\partial \varphi}{\partial x}
$$
on $\Delta^+ \cup ]-1,1[$.
Hence $\partial \varphi/\partial y$
is of class $\mathcal C^{1+1/2}$ on $\Delta^+ \cup ]-1,1[$, meaning that
$\varphi$ is of class $\mathcal C^{2+1/2}$ on $\Delta^+ \cup ]-1,1[$.
We prove now that $\varphi$ is of class $\mathcal C^{3+1/2}$ on
$\Delta^+ \cup ]-1,1[$. The reader will conclude, repeating the same
argument that $\varphi$ is of class $\mathcal C^\infty$ on
$\Delta^+ \cup ]-1,1[$.

Replace now the data $(M,J)$ and $\varphi$ by $(TM,J^c)$ and
$\varphi^c$ in Step three. The map $^2\varphi^c$ defined on $\Delta^+$ by
$^2\varphi^c(\zeta) = (\varphi^c(\zeta), d\varphi^c(\zeta)(v))$
is $^2J^c$-holomorphic on $\Delta^+$ ($^2J^c$ is the complete lift of $J^c$
to the second tangent bundle $T(TM)$. According to Step two, its
first derivative $\partial (^2\varphi^c)/\partial x$ is of class $C^{1+1/2}$ on
$\Delta^+ \cup ]-1,1[$. This means that the second derivatives
$\displaystyle \frac{\partial^2 \varphi}{\partial x^2}$ and
$\displaystyle \frac{\partial^2 \varphi}{\partial x \partial y}$
are $C^{1+1/2}$ on $\Delta^+ \cup ]-1,1[$. Differentiating
equation~(\ref{holo}) with respect to $y$, we prove that
$\displaystyle \frac{\partial^2 \varphi}{\partial y^2}$ is $C^{1+1/2}$ on
$\Delta^+ \cup ]-1,1[$ and so that $\varphi$ is $C^{3+1/2}$ on
$\Delta^+ \cup ]-1,1[$. \qed

\subsection{Boundary regularity of diffeomorphisms in wedges}
Let $\Gamma$ and $\Gamma'$ be two totally real maximal submanifolds in almost
complex manifodls $(M,J)$ and $(M',J')$. Let $W(\Gamma,M)$ be a wedge in
$M$ with edge $\Gamma$. 
\begin{proposition}\label{wed-reg}
If $F :W(\Gamma,M) \rightarrow M'$ is $(J,J')$-holomorphic and if the
cluster set of $\Gamma$ is contained in $\Gamma'$ then $F$ extends as a
$\mathcal C^\infty$ map up to $\Gamma$.
\end{proposition}

\noindent{\bf Proof of Proposition~\ref{wed-reg}.} In view of
Proposition~\ref{reflection} the proof is classical (see \cite{co-ga-su}). \qed

\vskip 0,1cm
As a direct application of Proposition~\ref{wed-reg} we obtain the following
partial version of Fefferman's Theorem :
\begin{corollary}\label{feff1}
Let $D$ and $D'$ be two smooth relatively compact domains in real manifolds.
Assume that $D$ admits an almost complex structure $J$ smooth on $\bar D$ and
 such that $(D,J)$ is strictly pseudoconvex. Let $f$ be a smooth
diffeomorphism $f:  D  \rightarrow  D'$, extending as a $\mathcal C^1$
diffeomorphism (still called $f$) between $\bar{D}$ and $\bar{D}'$.
Then $f$ is a smooth $\mathcal C^\infty$
diffeomorphism between $\bar D$ and $\bar D'$  if and only
if the direct image $f_*(J)$ of $J$ under
$f$ extends smoothly on $ \bar D'$ and $(D', f_*(J))$ 
is strictly pseudoconvex.
\end{corollary}

\noindent{\bf Proof of Corollary~\ref{feff1}.}
The cotangent lift $f^*$
of $f$ to the cotangent bundle over $D$, locally defined by
$f^*:=(f,^t(df)^{-1})$, is a $(\tilde{J},\tilde{J}')$-biholomorphism
from $T^*D$ to $T^*D'$, where $J':=f_*(J)$.
According to Proposition~\ref{prop-tot-real}, the conormal
bundle $\Sigma(\partial D)$ (resp. $\Sigma(\partial D')$) is a totally real
submanifold in $T^*M$ (resp. $T^*M'$).
%(see, for instance, L.Nirenberg-S.Webster-P.Yang \cite{nir-we-ya},
%S.Pinchuk-S.Khasanov \cite{pi-kh}, B.Coupet \cite{co} and
%F.Forstneric \cite{for}) and its main steps are explained in \cite{co-ga-su}.
We consider $\Sigma(\partial D)$ as the edge of a wedge
$W(\Sigma(\partial D),M)$ contained in $TD$. Then we may apply
Proposition~\ref{wed-reg} to $F=f^*$ to conclude. \qed

\section{Boundary estimates and the scaling process}
Our further considerations rely deeply on the following estimates
of the Kobayashi-Royden infinitesimal pseudometric obtained in~\cite{ga-su}~:
\vskip 0,1cm
\noindent{\bf Proposition A.}
{\it Let $D=\{\rho<0\}$ be a relatively compact domain in an almost complex
manifold $(M,J)$, where $\rho$ is a $\CC^2$ defining function of $D$,
strictly $J$-plurisubharmonic in a neighborhood of $\bar{D}$. Then
there exist positive constants $c$ and $\lambda_0$ such that for every
almost complex structure $J'$ defined in a neighborhood of $\bar{D}$ and 
such that $\|J'-J\|_{\CC^2(\bar{D})} \leq \lambda_0$ we have~:

\begin{equation}\label{e3}
K_{(D,J')}(p,v) \geq c\left[\frac{|\partial_J\rho(p)(v - iJ(p)v)|^2}
{|\rho(p)|^2} + 
\frac{\|v\|^2}{|\rho(p)|}\right]^{1/2},
\end{equation}
for every $p \in D$ and every $v \in T_pM$.}

\vskip 0,1cm
Let $D$ (resp. $D'$) be a strictly pseudoconvex domain in an almost complex
manifold $(M,J)$ (resp. $(M',J')$) and let $f$ be a $(J,J')$-biholomorphism
from $D$ to $D'$. Fix a point $p \in \partial D$ and a sequence $(p^k)_k$
in $D$ converging to $p$. After extraction we may assume that the sequence
$(f(p^k))_k$ converges to a point $p'$ in $\partial D'$. According to the
Hopf lemma, $f$ has the boundary distance property. Namely, there is a
positive constant $C$ such that
\begin{equation}\label{bdp}
(1/A) \ dist(f(p^k), \partial D') \leq dist(p^k, \partial D) \leq
A \ dist(f(p^k), \partial D'),
\end{equation}
where $A$ is independent of $k$ (see~\cite{co-ga-su}).

Since all our considerations are local we set $p=p'=0 \in \C^n$.
We may assume that $J(0) = J_{st}$ and $J'(0) = J_{st}$.
Let $U$  (resp. $V$) be a neighborhood of the origin in $\C^n$ such that
$D \cap U = \{z \in U : \rho(z,\bar{z}) : =
z_n + \bar{z}_n + Re(K(z)) + H(z,\bar{z}) + \cdots < 0\}$ 
(resp. $D' \cap V = \{w \in V : \rho'(w,\bar{w}) : =
w_n + \bar{w}_n + Re(K'(w)) + H'(w,\bar{w}) + \cdots < 0\}$) 
where 
$K(z)  = \sum k_{\nu\mu} z^{\nu}{z}^{\mu}$, $k_{\nu\mu} =
k_{\mu\nu}$,
$H(z) = \sum h_{\nu\mu} z^{\nu}\bar z^{\mu}$, $h_{\nu\mu} =
\bar h_{\mu\nu}$ and $\rho$ is a strictly $J$-plurisubharmonic function on
$U$ (resp. $K'(z)  = \sum k'_{\nu\mu} w^{\nu}{w}^{\mu}$, $k'_{\nu\mu} =
k'_{\mu\nu}$,
$H'(w) = \sum h'_{\nu\mu} w^{\nu}\bar w^{\mu}$, $h'_{\nu\mu} =
\bar h'_{\mu\nu}$ and $\rho'$ is a strictly $J'$-plurisubharmonic function on
$V$).

\subsection{Asymptotic behaviour of the tangent map of $f$}
We wish to understand the limit behaviour (when $k \rightarrow \infty$) of
$df(p^k)$. Consider the vector fields
$$
v^j:=(\partial \rho/\partial x^n)\partial / \partial x^j
- (\partial \rho/\partial x^j)\partial / \partial x^n
$$
for $j=1,\dots,n-1$, and
$$
v^n:=(\partial \rho/\partial x^n)\partial /
\partial y^n - (\partial \rho/\partial y^n)\partial / \partial x^n.
$$
Restricting $U$ if necessary, the vector fields $X^1,\dots,X^{n-1}$
defined by $X^j:=v^j-iJ(v^j)$ form a basis 
of the $J$-holomorphic tangent space to $\{\rho = \rho(z)\}$
at any $z \in U$. Moreover, if $X^n:=v^n-iJv^n$ then the family
$X:=(X^1,\dots,X^n)$ forms a basis of $(1,0)$ vector fields on $U$.
Similarly we define a basis
$X':=(X'^1,\dots,X'^n)$ of $(1,0)$ vector fields on $V$ such that
$(X'^1(w),\dots,X'^{n-1}(w))$ defines a basis
of the $J'$-holomorphic tangent space to $\{\rho' = \rho'(w)\}$
at any $w \in V$. 
We denote by $A(p^k):=(A(p^k)_{j,l})_{1 \leq j,l \leq n}$ the matrix of the map
$df(p^k)$ in the basis $X(p^k)$ and $X(f(p^k))$.

\begin{remark}\label{precise}
In sake of completeness we should write $X_0$ and $X'_0$ to emphasize that
the structure was normalized by the condition $J(0) = J_{st}$ and
$A(0,p^k)$ for $A(p^k)$. The same construction is valid for any
boundary point of $D$.
The corresponding notations will be used in Proposition~\ref{reality}.
\end{remark}
\begin{proposition}\label{tangent}
The matrix $A(p^k)$ satisfies the following estimates~:
$$
A(p^k)=\left(
\begin{array}{ccc}
O_{n-1,n-1}(1) & & O_{n-1,1}(dist(p^k,\partial D)^{-1/2})\\
 & & \\
O_{1,n-1}(dist(p^k,\partial D)^{1/2}) & & O_{1,1}(1)
\end{array}
\right).
$$
\end{proposition}

The matrix notation means that the following estimates are satisfied~:
$A(p^k)_{j,l} = O(1)$ for $1 \leq j,l \leq n-1$,
$A(p^k)_{j,n} = O(dist(p^k,\partial D)^{-1/2})$ for $1 \leq j \leq n-1$,
$A(p^k)_{n,l} = O(dist(p^k,\partial D)^{1/2})$ for $1 \leq l \leq n-1$ and
$A(p^k)_{n,n} = O(1)$. 
\vskip 0,1cm

The proof of Proposition~\ref{tangent} is given in \cite{co-ga-su}
(Proposition 3.5) in dimension two but is valid without any
modification in any dimension. This is based on Proposition~A. We note
that the asymptotic behaviour of $A(p^k)$ depends only on the distance
from the point to $\partial D$, not on the choice of the sequence
$(p^k)_k$.

\subsection{Scaling process and model domains}

The following construction is similar to the two dimensional case.
For every $k$ denote by $q^k$ the projection of $p^k$ to $\partial D$ and
consider the change of variables $\alpha^k$ defined by 
$$
\left\{
\begin{array}{ccccc}
(z^j)^* & = & \displaystyle \frac{\partial \rho}{\partial \bar z^n}(q^k)
(z^j - (q^k)^j)
- \displaystyle \frac{\partial \rho}{\partial \bar z^j}(q^k)(z^n - (q^k)^n),
& &
{\rm for} \ 1 \leq j \leq n-1,\\
(z^n)^* & = & \sum_{j=1}^n \displaystyle \frac{\partial \rho}{\partial z^j}
(q^k)(z^j - (q^k)^j).
\end{array}
\right.
$$
If $\delta_k := dist(p^k,\partial D)$ then $\alpha^k(p^k) = (0,-\delta_k)$
and $\alpha^k(D)=\{2Re z^n + O(\vert z
\vert^2) < 0\}$ near the origin. Moreover, the sequence $(\alpha^k)_*(J)$
converges to $J$ as $k \rightarrow \infty$, since the sequence
$(\alpha^k)_k$ converges
to the identity map. Let $(L^k)_k$ be a sequence of linear automorphisms of
$\R^{2n}$
such that $(T^k: = L^k
\circ \alpha^k)_k$ converges to the
identity, and $D^k:= T^k(D)$ is defined near the origin by
$D^k=\{\rho_k(z) = Re z^n + O(\vert z \vert^2) < 0\}$.
The sequence of almost complex structures
$(J_k:= (T^k)_*(J))_k$ converges to $J$ as $k \rightarrow \infty$
and $J_k(0) = J_{st}$.
Furthermore $\tilde p_k := T^k(p^k)$ satisfies
$\tilde p_k = (o(\delta_k),\delta_k'' + io(\delta_k))$
with
$\delta_k''\sim \delta_k$. 
%We denote by $\Gamma^k = \{\rho_k = 0 \}$ the
%image of $\partial D$ under $T^k$.

We proceed similarly on $D'$. 
We denote by $s^k$ the projection of $f(p^k)$ onto $\partial D'$ and
we define the transformation $\beta^k$ by
$$
\left\{
\begin{array}{ccccc}
(w^j)^* & = & \displaystyle \frac{\partial \rho'}{\partial \bar w^n}(s^k)
(w^j - (s^k)^j)
- \displaystyle \frac{\partial \rho'}{\partial \bar w^j}(s^k)(w^n - (s^k)^n),
& & {\rm for} \ 1 \leq j \leq n-1,\\
(w^n)^* & = & \sum_{j=1}^n \displaystyle \frac{\partial \rho'}{\partial w^j}
(s^k)(w^j - (s^k)^j).
\end{array}
\right.
$$
We define a sequence $(T'^k)_k$ of
linear transformations converging to the identity and satisfying the
following properties. The domain 
$(D^k)':= T'^k(D')$ is defined near the origin by
$(D^k)'=\{\rho_k'(w) := Re w^n + O(\vert w \vert^2) < 0\}$,
and $\tilde f(p_k) = T'^k(f(p_k)) =
(o(\varepsilon_k),\varepsilon_k''+ io(\varepsilon_k))$
with $\varepsilon_k'' \sim \varepsilon_k$, where
$\varepsilon_k = dist(f(p_k),\partial D')$. 
The sequence of almost complex structures $(J_k':= (T'^k)_*(J'))_k$
converges to $J'$ as $k \rightarrow \infty$ and $J_k'(0) = J_{st}$.  

Finally, the map $ f^k:= T'^k \circ f \circ (T^k)^{-1}$ satisfies
$f^k(\tilde p_k) = \tilde f(p_k)$  and is
a $(J_k,J'_k)$-biholomorphism between the domains $D^k$ and $(D')^k$. 

Let $\phi_k : ('z,z^n) \mapsto (\delta_k^{1/2} {'z},\delta_kz^n)$ and
$\psi_k(w',w^n)=
(\varepsilon_k^{1/2}w',\varepsilon_kw^n)$ and set $\hat f^k =
(\psi_k)^{-1} \circ f^k \circ \phi_k$.
The map $\hat f^k$ is $(\hat J_k,\hat J'_k)$-biholomorphic, where
$\hat J_k:=((\phi_k)^{-1})_*(J_k)$ and
$\hat J'_k:= (\psi_k^{-1})_*(J'_k)$.
If $\hat D^k:=\phi_k^{-1}(D^k)$ and
$(\hat{D'})^k:=\psi_k^{-1}((D')^k)$ then
$\hat D^k = \{ z \in \phi_k^{-1}(U): \hat \rho_k(z) < 0\}$
where
$$
\hat \rho_k(z): = \delta_k^{-1}\rho(\phi_k(z)) = 2Re z^n + \delta_k^{-1}[2
Re K(\delta_k^{1/2}{'z},\delta_kz^n) + H(\delta_k^{1/2}{'z},\delta_kz^n)
+  o(\vert (\delta_k^{1/2}{'z},\delta_kz^n) \vert^2).
$$
and $(\hat D')^k=\{w \in \phi_k^{-1}(V): \hat \rho'_k(z) < 0\}$
where
$$
\hat \rho'_k(w): = \varepsilon_k^{-1}\rho'(\psi_k(w)) = 2Re w^n +
\varepsilon_k^{-1}[2 Re K'(\varepsilon_k^{1/2}w',\varepsilon_kw^n) +
H'(\varepsilon_k^{1/2}w',\varepsilon_kw^n)
+  o(\vert (\varepsilon_k^{1/2}w',\varepsilon_kw^n) \vert^2).
$$
Since $U$
is a neighborhood of the origin, the pullbacks $\phi_k^{-1}(U)$
converge to $\C^n$ and the functions $\hat\rho_k$ converge
to $\hat \rho(z) = 2Re z^n + 2Re K({'z},0) + H({'z},0)$ in the $\CC^2$ norm
on compact subsets of $\C^n$. Similarly, since $V$
is a neighborhood of the origin, the pullbacks $\psi_k^{-1}(U')$
converge to $\C^n$ and the functions $\hat\rho_k'$ converge
to $\hat \rho'(w) = 2Re w^n + 2Re K'({'z},0) + H'({'z},0)$ in the $\CC^2$ norm
on compact subsets of $\C^n$. If $\Sigma :=
\{z \in \C^n: \hat \rho(z) < 0 \}$ and $\Sigma' := \{w \in \C^n:
\hat \rho'(w) < 0 \}$ the sequence of points $\hat p_k =
\phi_k^{-1}(\tilde p_k) \in \hat D^k$ converges to the point $(0,-1) \in
\Sigma$ and the sequence of points $\hat f(p_k) =
\psi^{-1}_k(\tilde f(p_k)) \in \hat{D'}^k$ converges to $(0,-1) \in
\Sigma'$. Finally $\hat{f}^k(\hat p_k) = \hat f(p_k)$.

The limit behaviour of the dilated objects is given by the following
proposition.
\begin{proposition}\label{convseq}
$(i)$ The sequences $(\hat J_k)$ and $(\hat J'_k)$ of almost complex
structures converge to model structures $J_0$ and $J'_0$
uniformly (with all partial derivatives of any order) on compact subsets of
$\C^n$.

\vskip 0,1cm
$(ii)$ $(\Sigma,J_0)$ and $(\Sigma',J'_0)$ are model domains.

\vskip 0,1cm
$(iii)$ The sequence $(\hat f^k)$ (together with all derivatives) is a
relatively compact family (with respect to the compact open topology) on
$\Sigma$; every cluster point $\hat f$ is
a $(J_0,J'_0)$-biholomorphism between $\Sigma$ 
and $\Sigma'$, satisfying $\hat f(0,-1) = (0,-1)$ and
%$(\partial \hat f^n/\partial z^n)(0,-1) = 1$.
$\hat f^n('0,z^n) = z^n$ on $\Sigma$.
\end{proposition}

\noindent{\it Proof of Proposition~\ref{convseq}.}

\vskip 0,1cm
Proof of $(i)$. We focus on structures $\hat{J}_k$.
Consider $J=J_{st} + L(z) +
  O(|z|^2)$ as a matrix
valued function, where $L$ is a real linear matrix.
The Taylor expansion of $J_k$ at
the origin is given by  $J_k = J_{st} + L^k(z) +   O(|z|^2)$
on $U$, uniformly with respect to $k$. Here $L^k$ is a real linear
matrix converging to $L$ at infinity. 
Write $\hat{J}_k = J_{st} + \hat{L}^k +   O(\delta_k)$.
If $L^k=(L^k_{j,l})_{j,l}$ then
$\hat{L}^k_{j,l}= L^k_{j,l}(\phi_k(z))$ for $1 \leq j \leq n-1,\ 1 \leq
l \leq n$, $\hat{L}^k_{n,l}=\delta_k^{-1/2}L^k{n,l}(\phi_k(z))$
for $1 \leq l \leq n-1$ and $\hat{L}^k_{n,n}=L^k_{n,n}(\phi_k(z))$.
This gives the conclusion. 

\vskip 0,1cm
Proof of $(ii)$. We focus on $(\Sigma,J_0)$.
By the invariance of the Levi form we
have ${\mathcal L}^{J_k}(\rho_k)(0)(\phi_k(v)) = {\mathcal
L}^{\hat J_k}(\rho_k \circ \phi_k)(0)(v)$.
Write $J_0 = J_{st} + L^\infty$.
Since $\rho_k$ is strictly $J_k$-plurisubharmonic uniformly with respect to $k$
($\rho_k$ converges to $\rho$ and $J_k$ converges to $J$),
multiplying by $\delta_k^{-1}$ and
passing to the limit at the right side as $k \rightarrow \infty$, 
we obtain that
${\mathcal L}^{J_0}(\hat \rho)(0)(v) \geq 0$ for any $v$. Now let $v =
(v',0) \in T_0(\partial \Sigma)$. Then
$\phi_k(v) = \delta_k^{1/2}v$ and so 
${\mathcal L}^J_k(\rho)(0)(v) = {\mathcal
  L}^{\hat J_k}(\rho_k)(0)(v)$. Passing to the limit as $k$
tends to infinity, we obtain that
${\mathcal L}^{J_0}(\hat \rho)(0)(v) > 0$
for any $v = (v',0)$ with $v' \neq 0$. 

\vskip 0,1cm
Proof of $(iii)$. This statement is a consequence of Proposition~A.
We refer to Section~7 of \cite{co-ga-su} for the existence
and the biholomorphy of $\hat{f}$. We prove the identity on $\hat f^n$.
Let $t$ be a real positive number. Then we have~:
\begin{lemma}\label{infty}
$\lim_{t \rightarrow \infty} \hat{\rho}'(\hat{f}('0,-t)) =  \infty$.
\end{lemma}
\noindent{\it Proof of Lemma~\ref{infty}}.
According to the boundary distance property~(\ref{bdp}) we have
$$
|\rho'(f \circ (T^k)^{-1} \circ \phi_k)('0,-t)| \geq C \ 
dist(T_k^{-1}('0,-\delta_k t)).
$$
Then
$$
|\hat{\rho}'_k(\hat{f}^k('0,-t))| \geq C \varepsilon_k^{-1}\delta_k \ t.
$$
Since $\hat{\rho}'_k$ converges to $\hat{\rho}'$ uniformly on
compact subsets of $\Sigma'$ and $\varepsilon_k \simeq \delta_k$ (by
the boundary distance property~(\ref{bdp})) we obtain~:
$$
|\hat{\rho}'(\hat{f}('0,-t))| \geq Ct.
$$
This proves Lemma~\ref{infty}. \qed

\vskip 0,1cm
We turn back to the proof of part $(iii)$ of Proposition~\ref{convseq}.
Assume first that $J$ (and similarly $J'$) are not integrable
(see Proposition~\ref{prop-hyp}). Consider a $J$-complex hypersurface
$A \times \C$ in $\C^n$ where $A$ is a $J_{st}$ complex hypersurface in
$\C^{n-1}$.
Since $f((A \times \C) \cap \mathbb H_{P_1}) = (A' \times \C) \cap
\mathbb H_{P_2}$ where $A'$ is a $J_{st}$ complex hypersurface in
$\C^{n-1}$, it follows that the restriction of $\hat f^n$ to $\{'z='0,
Re(z^n) < 0\}$ is a $J_{st}$ automorphism of $\{'z='0, Re(z^n) < 0\}$.
Let $\phi : \zeta \mapsto (\zeta -1)/(\zeta + 1)$. The function
$\hat{g}:=\phi^{-1} \circ \hat{f}^n \circ \phi$ is a $J_{st}$ automorphism
of the unit disc in $\C$. In view of Lemma~\ref{infty} this satisfies
$\hat{g}(0) = 0$ and $\hat{g}(1) = 1$. Hence $\hat{g} \equiv id$ and
$\hat{f}^n('0,z^n) = z^n$ on $\Sigma$.

Assume now that $J$ and $J'$ are integrable.
Let $F$ (resp.
$F'$) be the diffeomorphism from $\Sigma$ to $\mathbb H_{P}$
(resp. from $\Sigma$ to $\mathbb H_{P'}$) given in the proof of
Proposition~\ref{prop-hyp}. The diffeomorphism $g:=F' \circ f \circ
F^{-1}$ is a $J_{st}$-biholomorphism from $\mathbb H_{P}$ to $\mathbb
H_{P'}$ satisfying $g('0,-1) = ('0,-1)$. Since $(\Sigma,J)$ and
$(\Sigma', J')$ are model domains, the domains $\mathbb H_{P}$ and
$\mathbb H_{P'}$ are strictly $J_{st}$-pseudoconvex. In particular, since
$P$ and $P'$ are homogeneous of degree two, there are linear complex maps
$L,\ L'$ in $\C^{n-1}$ such that the map $G$ (resp. $G'$) defined by
$G('z,z_n)=(L('z),z_n)$ (resp. $G'('z,z_n)=(L'('z),z_n)$) is a
biholomorphism from $\mathbb H_{P}$ (resp. $\mathbb H_{P'}$) to
$\mathbb H$. The map $G' \circ g \circ G^{-1}$ is an automorphism of
$\mathbb H$ satisfying $G' \circ g \circ G^{-1}('0,-1) = ('0,-1)$.
Let $\Phi$ be the $J_{st}$ biholomorphism from $\mathbb H$ to
the unit ball $\mathbb B_n$ of $\C^n$ defined by
$\Phi('z,z^n) = (\sqrt{2}'z/1-z^n,(1+z^n)/(1-z^n))$.
Let $\hat{g} :=\Phi^{-1} \circ g
\circ \Phi$. In view of lemma~\ref{bdp} this satisfies
$\hat{g}(0) = 0$ and $\hat{g}('0,1)=('0,1)$. Hence $\hat{g}^n \equiv id$
and $\hat{f}^n('z,z^n) = z^n$ for every $z$ in $\Sigma$. \qed

\vskip 0,1cm
According to part $(ii)$ of Proposition~\ref{convseq}
and restricting $U$ if necessary, one may view
$D \cap U$ as a strictly $J_0$-pseudoconvex domain
in $\C^n$ and $J$ as a small deformation of ${J}_0$ in a
neighborhood of $\bar{D} \cap U$. The same holds for $D' \cap V$.

\vskip 0,1cm
For $p \in \partial D$ and $z \in D$ let $X_p(z)$ and $X'_{f(p)}(f(z))$
be the basis of $(1,0)$ vector fields defined in Subsection~3.2
(see Remark~\ref{precise}).
The elements of the matrix of $df_z$ in
the bases $X_p(z)$ and $X'_{f(p)}(f(z))$ are denoted by
$A_{js}(p,z)$. According to Proposition~\ref{tangent} the function
$A_{n,n}(p,\cdot)$ is upper bounded on $D$. 

\begin{proposition}
\label{reality}
We have:
\begin{itemize}
\item[(a)] Every cluster point of the function $z \mapsto A_{n,n}(p,z)$
is real when $z$ tends to $p \in \partial D$.
\item[(b)] For $z \in D$, let $p \in \partial D$ such that
$|z-p| = dist(z,\partial D)$. There exists a constant $A$, independent of
$z \in D$, such that $\vert A_{n,n}(p,z) \vert \geq A$.
\end{itemize}
\end{proposition}

\vskip 0,1cm
\noindent{\it Proof of Proposition~\ref{reality}}.
(a) Suppose that there exists a sequence of points $(p^k)$ converging
to a boundary point $p$ such that $A_{n,n}(p,\cdot)$ tends to a complex number
$a$. Applying the above scaling construction,
we obtain a sequence of maps $(\hat f^k)_k$.
For $k \geq 0$ consider the dilated vector fields
$$
Y^j_k:=\delta_k^{1/2}((\phi_k^{-1}) \circ T^k)(X^j(p^k))
$$
for $j=1,\dots,n-1$, and
$$
Y^n_k:=\delta_k((\phi_k^{-1})\circ T^k)(X_n(p^k)).
$$
Similarly we define
$$
Y'^j_k:=\varepsilon_k^{-1/2}((\psi_k^{-1}) \circ T'^k)
(X'^j(f(p^k)))
$$
for $j=1,\dots,n-1$, and
$$
Y'^n_k:=\varepsilon_k^{-1}((\psi_k^{-1})\circ T'^k)(X'_n(f(p^k))).
$$
For every $k$, the $n$-tuple
$Y^k:= (Y^1_k,\dots,Y^n_k)$ is a basis of $(1,0)$ vector fields for
the dilated structure $\hat{J}^k$. In view of Proposition~\ref{convseq}
the sequence $(Y^k)_k$
converges to a basis
of $(1,0)$ vector fields of $\C^n$ (with respect to $J_0$) as $k$
tends to $\infty$. Similarly, the $n$-tuple
$Y'^k := (Y'^1_k,\dots,Y'^n_k)$ is a basis of $(1,0)$ vector fields for
the dilated structure $\hat{J}'^k$ and $(Y'^k)_k$
converges to a basis of $(1,0)$ vector fields
of $\C^n$ (with respect to $J'_0$) as $k$ tends to $\infty$.
In particular the last components $Y^n_k$ and $Y'^n_k$
converge to the $(1,0)$ vector field $\partial / \partial z^n$.
Denote by $\hat A^k_{js}$ the elements of the
matrix of $d\hat f^k(0,-1)$. Then $A^k_{n,n}$ converges to $(\partial
\hat f^n/\partial z^n)(0,-1) = 1$, according to Proposition~\ref{convseq}.
On the other hand, $A^k_{n,n} = \varepsilon_k^{-1}\delta_k A_{n,n}$
converges to $a$ by the boundary distance preserving property~(\ref{bdp}).
This gives the statement.

(b) Suppose that there is a sequence of points $(p^k)$ converging
to the boundary such that $A_{n,n}$ tends to $0$. Repeating precisely
the argument of (a), we obtain that $(\partial \hat f^n/\partial
z^n)(0,-1) = 0$; this contradicts part $(iii)$ of Proposition~\ref{convseq}.
\qed

\section{Proof of Theorem~\ref{theo-fefferman}}
From now on we are in the hypothesis of Theorem~\ref{theo-fefferman}.
The key point of the proof of Theorem~\ref{theo-fefferman} consists in the
following claim :

\vskip 0,1cm
{\it Claim : The cluster set of the cotangent lift $f^*$ on
$\Sigma(\partial D)$ is contained in $\Sigma(\partial D')$.}

\vskip 0,1cm
Indeed, assume for the moment the claim satisfied.
We recall that according to Proposition~\ref{prop-tot-real}
the conormal bundle $\Sigma_J(\partial D)$ of $\partial D$ is a
totally real submanifold in the cotangent bundle $T^*M$.
Consider the set $S = \{(z,L) \in \R^{2n} \times \R^{2n} :
dist((z,L),\Sigma_J(\partial D)) \leq dist(z,\partial D), z \in D \}$.
Then, in a neighborhood $U$ of any totally real point of
$\Sigma_J(\partial D)$, the set S contains a wedge $W_U$ with
$\Sigma_J(\partial D) \cap U$ as totally real edge.

Then in view of Proposition~\ref{wed-reg} we obtain the following
Proposition~:
\begin{proposition}
\label{wedges}
Restricting the aperture of the wedge $W_U$ if necessary,
the map $f^*$ extends to $W_U \cup \Sigma(\partial D)$
as a $\CC^{\infty}$-map. 
\end{proposition}
Proposition~\ref{wedges} implies immediately that $f$ extends
as a smooth $\mathcal C^\infty$ diffeomorphism from $\bar{D}$ to
$\bar{D'}$.

Therefore the proof of Theorem~\ref{theo-fefferman} can be reduced to the
proof of the claim.

\vskip 0,1cm
\noindent{\it Step one.}
We first reduce the problem to the following
local situation. Let $D$ and $D'$ be domains in $\C^n$, $\Gamma$ and
$\Gamma'$ be open $\CC^{\infty}$-smooth pieces of their boundaries,
containing the origin. We assume that an almost complex structure $J$
is defined and $\CC^{\infty}$-smooth in a neighborhood of the closure
$\bar D$, $J(0) = J_{st}$.
Similarly, we assume that $J'(0) = J_{st}$. The hypersurface
$\Gamma$ (resp. $\Gamma'$) is supposed to be strictly $J$-pseudoconvex
(resp. strictly $J'$-pseudoconvex). Finally, we assume that $f: D
\rightarrow D'$ is a $(J,J')$-biholomorphic map. It follows from the estimates
of the Kobayashi-Royden infinitesimal pseudometric given in \cite{ga-su}
that $f$ extends as a $1/2$-H\"older
homeomorphism between  $D \cup \Gamma$ and $D' \cup \Gamma'$, such that
$f(\Gamma) = \Gamma'$ and $f(0) = 0$. Finally
$\Gamma$ is defined in a neighborhood of the origin
by the equation $\rho(z) = 0$ where $\rho(z) = 2Re
z^n + 2Re K(z) + H(z) + o(\vert z \vert^2)$ and $K(z) = \sum
K_{\mu\nu}z^{\mu\nu}$, $H(z) = \sum h_{\mu\nu}z^{\mu}\bar
z^{\nu}$, $k_{\mu\nu} = k_{\nu\mu}$, $h_{\mu\nu} = \bar
h_{\nu\mu}$. As we noticed at the end of Section~3 the hypersurface
$\Gamma$ is strictly $\hat{J}$-pseudoconvex at the origin. The hypersurface
$\Gamma'$ admits a similar local representation. In what follows we
assume that we are in this setting. 

Let $\Sigma := \{ z \in \C^n: 2Re z^n + 2Re K('z,0) + H('z,0) < 0\}$,
$\Sigma' := \{ z \in \C^n: 2Re z^n + 2Re K'('z,0) + H'('z,0) < 0\}$.
If $(p^k)$ is a sequence of points in $D$  converging to $0$, then according
to Proposition~\ref{convseq}, the
scaling procedure associates with the pair $(f,(p^k)_k)$ two linear almost
complex structures ${J}_0$ and ${J}'_0$, both defined on $\C^n$,
and a $(J_0,{J}'_0)$-biholomorphism $\hat{f}$ between $\Sigma$ and
$\Sigma'$. Moreover $(\Sigma,J_0)$ and $(\Sigma',J'_0)$ are model 
domains. To prove that the cluster set of the cotangent lift of $f$ at a point
in $N(\Gamma)$ is contained in $N(\Gamma')$, it is sufficient to prove that
$(\partial \hat{f}^n / \partial z^n)('0,-1) \in \mathbb R \backslash \{0\}$.

\vskip 0,1cm
\noindent{\it Step two.} The proof of the Claim is given by the following
Proposition.

\begin{proposition}
\label{cluster}
Let $K$ be a compact subset of the totally real part of the conormal
bundle $\Sigma_J(\partial D)$. Then the cluster set of the cotangent lift
$f^*$ of $f$ on the conormal bundle 
$\Sigma(\partial D)$, when $(z,L)$ tends to $\Sigma_J(\partial D)$
along the wedge $W_U$, is relatively compactly contained 
in the totally real part of $\Sigma(\partial D')$.
\end{proposition}
We recall that the totally real part of $\Sigma(\partial D')$ is
the complement of the zero section in $\Sigma(\partial D')$.

\vskip 0,1cm
\noindent{\sl Proof of Proposition~\ref{cluster}}.  Let $(z^k,L^k)$ be
a sequence in $W_U$ converging to $(0,\partial_J\rho(0)) =
(0,dz^n)$.  We shall prove that the sequence of
linear forms $Q^k := {}^tdf^{-1}(w^k)L^k$, where $w^k = f(z^k)$, converges
to a linear form which up to a {\it real} factor (in view of Part (a)
of Proposition \ref{reality}) coincides with $\partial_{J} \rho(0)=
dz^n$ (we recall that ${}^t$ denotes the transposed map).  It is
sufficient to prove that the $(n-1)$ first component of $Q^k$ with respect to
the dual basis $(\omega_1,\dots,\omega_n)$ of $X$ converge to $0$ and the
last one is bounded below from the origin as $k$ goes to infinity.
The map $X$ being of class $\CC^1$ we can replace $X(0)$ by $X(w^k)$.
Since $(z^k,L^k) \in W_U$, we have $L^k = \omega_n(z^k) +
O(\delta_k)$, where $\delta_k$ is the distance from $z^k$ to the
boundary. Since $\vert\vert\vert df^{-1}_{w^k} \vert\vert\vert =
0(\delta_k^{-1/2})$, we have $Q^k = {}^tdf^{-1}_{w^k}(\omega_n(z^k)) +
O(\delta_k^{1/2})$.  By Proposition~\ref{convseq}, the components of
${}^tdf^{-1}_{w^k}(\omega_n(z^k))$ with respect to the basis
$(\omega_1(z^k),\dots,\omega_n(z^k))$ are the elements of the last line of
the matrix $df^{-1}_{w^k}$ with respect to the basis $X'(w^k)$ and
$X(z^k)$. So its $(n-1)$ first components are $0(\delta_k^{1/2})$ and
converge to $0$ as $k$ tends to infinity. Finally the component $A_{n,n}^k$
is bounded below from the origin by Part (b) of
Proposition~\ref{reality}. \qed

\section{Compactness principle}
In this section we prove Theorem~\ref{wr}.

We note that condition $(ii)$ is equivalent to the existence,
at each $p \in \partial D$, of a strictly $J$-plurisubharmonic local
defining function for $\partial D$ (consider the function $\rho + C \rho^2$
for a sufficiently large positive $C$).

We first recall the following result proved in \cite{ga-su} :

\vskip 0,1cm
\noindent{\bf Proposition B}. {\it (Localization principle)
Let $D$ be a domain in an almost complex
manifold $(M,J)$, let $p \in \bar{D}$, let $U$ be a neighborhood of
$p$ in $M$ (not necessarily contained in $D$) and let $z:U \rightarrow
\mathbb B$ be the diffeomorphism given by Lemma~\ref{suplem1}.  
Let $u$ be a $\mathcal C^2$ function on $\bar{D}$, negative and
$J$-plurisubharmonic on $D$. We assume that $-L \leq u < 0$ on $D \cap
U$ and that $u-c|z|^2$ is $J$-plurisubharmonic on $D \cap U$, where
$c$ and $L$ are positive constants. Then there exist a positive
constant $ s$ and a neighborhood $V \subset \subset U$ of $p$,
depending on $c$ and $L$ only, such that for $q \in D \cap V$ and $v
\in T_qM$ we have the following estimate:

\begin{equation}\label{e2}
K_{(D,J)}(q,v) \geq s K_{(D \cap U,J)}(q,v).
\end{equation}
}

\vskip 0,1cm
We can now prove Theorem~\ref{wr}.

\vskip 0,1cm
\noindent{\bf Proof of Theorem~\ref{wr}}.
We assume that the assumptions of Theorem~\ref{wr} are
satisfied. We proceed by contradiction.
Assume that there is a compact $K_0$ in $M$, points $p^\nu \in M$ and a point
$q \in \partial D$ such that $\lim_{\nu \rightarrow \infty}f^\nu(p^\nu) = q$.

\begin{lemma}\label{met-kob}
For every relatively compact neighborhood $V$ of $q$ there is $\nu_0$
such that for $\nu \geq \nu_0$ we have~:
$\lim_{x \rightarrow q}inf_{q' \in D \cap \partial V}d^K_{(D,J_\nu)}=\infty$.
\end{lemma}

\noindent{\it Proof of Lemma~\ref{met-kob}}. Restricting $U$ if
necessary, we may assume that the function $\rho + C \rho ^2$ is a
strictly $J_\nu$-plurisubharmonic function in a neighborhood of
$\bar{D} \cap U$, for sufficiently large $\nu$.
Moreover, using Proposition~B, we can focus on $K_{D \cap U}$. Smoothing
$D \cap U$, we may assume that the hypothesis of Proposition~A are satisfied
on $D \cap U$, uniformly for sufficiently large $\nu$.
In particular, the inequality~(\ref{e3}) is satisfied on
$D \cap U$, with a positive constant $c$ independent of $\nu$.
The result follows by a direct integration of this inequality.
\qed
 
\vskip 0,1cm
The following Lemma is a corollary of Lemma~\ref{met-kob}.
\begin{lemma}\label{lem3.3.1}
For every $K \subset \subset M$ we have : 
$\lim_{\nu \rightarrow \infty}f^\nu(K) =q$.
\end{lemma}

\noindent{\it Proof of Lemma \ref{lem3.3.1}}. Let $K \subset \subset M$ 
be such that $x^0 \in
K$. Since the function $x \mapsto d_D^K(x^0,x)$ is bounded from above by a
constant $C$ on $K$, it follows from the decreasing property of the Kobayashi 
pseudodistance that

\begin{equation}\label{eq2}
d_{(D,J_\nu)}^K(f^\nu(x^0),f^\nu(x)) \leq C
\end{equation}
for every $\nu$ and every 
$x \in K$. It follows from Lemma~\ref{met-kob} that for
every $V \subset \subset U$, containing $p$, we have :
\begin{equation}\label{eq3}
\lim_{\nu \rightarrow \infty}d_{(D,J_\nu)}^K 
(f^\nu(x^0),D \cap \partial V) = +\infty.
\end{equation} 
Then from conditions (\ref{eq2}) and (\ref{eq3}) we deduce that 
$f^\nu(K) \subset V$ for every sufficiently large $\nu$. 
This gives the statement. \qed 

\vskip 0,1cm
Fix now a point $p \in M$ and denote by $p^\nu$ the point $f^\nu(p)$.
We may assume that the sequence $(J_\nu:=f^\nu_*(J))_\nu$ converges
to an almost complex structure $J'$ on $\bar{D}$ and according to
Lemma~\ref{lem3.3.1} we may assume that
$\lim_{\nu \rightarrow \infty}p^\nu = q$.
We apply Subsection~4.3 to the domain $D$ and the sequence $(q^\nu)_\nu$.
We denote by $T^\nu$ the linear transformation
$T^\nu:=M^\nu \circ L^\nu \circ \alpha^\nu$, as in Subsection~4.3, and
we consider $D^\nu:=T^\nu(D)$, and $J^\nu:=T^\nu_*(J_\nu)$.
If $\phi_\nu$ is the nonisotropic dilation $\phi_\nu:('z,z^n) \mapsto
(\delta_\nu^{1/2}\ 'z,\delta_\nu z^n)$ then we set
$\hat{f}^\nu:=\phi_\nu^{-1} \circ T^\nu \circ f$ and
$\hat{J}^\nu:=(\phi_\nu^{-1})_*(J^\nu)$. We also consider
$\hat{\rho}_\nu:=\delta_\nu^{-1} \circ \rho \circ \phi_\nu$
and $\hat{D}^\nu:=\{\hat{\rho}_\nu < 0\}$.
As proved in Subsection~4.3, the sequence $(\hat{D}^\nu)_\nu$ converges,
in the local Hausdorff convergence, to a domain
$\Sigma:=\{z \in C^n:\hat \rho(z) := 2Re z^n + 2Re K({'z},0) + H({'z},0)<0\}$,
where $K$ and $H$ are homogeneous of degree two.
According to Proposition~\ref{convseq} we have~:

$(i)$ The sequence $(\hat{J}^\nu)$ converges to a model almost complex
structure $J_0$, uniformly (with all partial derivatives of any order)
on compact subsets of $\C^n$,

$(ii)$ $(\Sigma,J_0)$ is a model domain,

$(iii)$ the sequence $(\hat{f}^\nu)_\nu$ converges to a $(J,J_0)$
holomorphic map $F$ from $M$ to $\Sigma$.

\vskip 0,1cm
To prove Theorem~\ref{wr}, it remains to prove that $F$ is a diffeomorphism
from $M$ to $\Sigma$. 
We first notice that according to condition $(ii)$ of Theorem~\ref{wr}
and Lemma~\ref{met-kob}, the domain $D$ is
complete $J_\nu$-hyperbolic. In particular, since $f^\nu$ is a $(J,J_\nu)$
biholomorphism from $M$ to $D$, the manifold $M$ is complete $J$-hyperbolic.
Consequently, for every compact subset $L$ of $M$, there is a positive
constant $C$ such that for every $z \in L$ and every $v \in T_zM$ we have
$K_{(M,J)}(z,v) \geq C\|v\|$.
Consider the map $\hat{g}^\nu:=(\hat{f}^\nu)^{-1}$.
This is a $(\hat{J}^\nu,J)$ biholomorphism from $\hat{D}^\nu$ to $M$.
Let $K$ be a compact set in $\Sigma$. We may consider $\hat{g}^\nu(K)$
for sufficiently large $\nu$. By the decreasing property of the Kobayashi
distance, there is a compact subset $L$ in $M$ such that
$\hat{g}^\nu(K) \subset L$ for sufficiently large $\nu$. Then for every
$w \in K$ and for every $v \in T_w\Sigma$ we obtain, by the decreasing of the
Kobayashi-Royden infinitesimal pseudometric~:
$$
\|df^\nu(w)(v)\| \leq (1/C) \|v\|,
$$
uniformly for sufficiently large $\nu$.
According to Ascoli Theorem, we may extract from
$(\hat{g}^\nu)_\nu$ a subsequence, converging to a map $G$ from
$\Sigma$ to $M$. Finally, on any compact subset $K$ of $M$, by
the equality $\hat{g}^\nu \circ \hat{f}^\nu = id$ we obtain $F \circ G = id$.
This gives the result. \qed

\vskip 0,1cm
As a corollary of Theorem~\ref{wr} we obtain the following almost complex
version of the Wong-Rosay Theorem in real dimension four~:
\begin{corollary}\label{wr-2}
Let $(M,J)$ (resp. $(M',J')$) be an almost complex manifold of real dimension
four. Let $D$ (resp. $D'$) be a relatively compact domain in $M$ (resp. $N$).
Consider a sequence $(f^\nu)_\nu$ of diffeomorphisms from $D$ to
$D'$ such that the sequence $(J_\nu:=f^\nu_*(J))_\nu$ extends to $\bar{D}'$
and converges to $J'$ in the $C^2$ convergence on $\bar{D}'$.

Assume that there is a point $p\in D$ and a point $q \in \partial D'$ such
that $\lim_{\nu \rightarrow \infty}f^\nu(p) = q$ and such that
$D'$ is strictly $J'$-pseudoconvex at $q$.
Then there is a $(J,J_{st})$-biholomorphism from $M$ to the unit ball $\B^2$
in $\C^2$. 
\end{corollary}

\noindent{\bf Proof of Corollary~\ref{wr-2}}.
The proof of Corollary~\ref{wr-2} follows exactly the same lines.
Indeed, by assumtion there is a fixed neighborood $U$ of $q$ such that
$D' \cap U$ is strictly $J_\nu$-pseudoconvex on $U$.
According to Lemma~\ref{lem3.3.1}, we know that for every compact subset
$K$ of $D$ the set $f^\nu(K)$ is contained in $V$ for sufficiently large
$\nu$. If we fix a point $p \in D$ we may therefore apply Subsection~4.4
to the sequence $(f^\nu(p))_\nu$ and to the domain $D'$ (with $V$ exactly as
in Subsection~4.4). The proof is then identical to the proof of
Theorem~\ref{wr}. \qed

\end{document}